\documentclass[draft,amsmath,amscd,amsbsy,amssymb]{amsart}
 \relax
\chardef\bslash=`\\ % p.  424, TeXbook %\newcommand{\ntt}{\seriesm\shape n\tt}

\makeatletter \def\verbatim{\interlinepenalty\@M \@verbatim
\leftskip\@totalleftmargin\advance\leftskip2pc \frenchspacing\@vobeyspaces
\@xverbatim} \makeatother \hfuzz1pc
\makeatletter \def\dgt@k{\dg@DX=-3 \dg@DY=2 \dg@SIZE=3} \makeatother
\makeatletter \def\dgt@kk{\dg@DX=3 \dg@DY=-1 \dg@SIZE=3}% \makeatother
\theoremstyle{plain}
\newtheorem{thm}{Theorem}[section]

\newtheorem{lemma}[thm]{Lemma}

\theoremstyle{definition}

\newcommand{\R}{{\mathbb R}}

\newcommand{\U}{{\mathcal U}}

\numberwithin{equation}{section}

\newcommand{\diam}{{\mathrm d}{\mathrm i} {\mathrm a}{\mathrm m}\ }
\newcommand{\asd}{{\mathrm a}{\mathrm s} {\mathrm d}{\mathrm i}{\mathrm m}}

\newcommand{\N}{{\mathbb N}}

\newcommand{\e}{{\varepsilon}}

\begin{document}
%%%%%%% Begin Topmatter %%%%%%%%%%
\title[Hyperspaces of probability measures and  extensors]
{Convex hyperspaces of probability measures and extensors
 in the asymptotic category}
\author[D.~Repov\v s]{Du\v san Repov\v s}
\address{Faculty of Mathematics and Physics, and Faculty of Education,
University of Ljubljana, P.O.B. 2964, Ljubljana, 1001, Slovenia}
\email{dusan.repovs@guest.arnes.si}

\author[M.~Zarichnyi]{Mykhailo Zarichnyi}
\address{Department of Mechanics and Mathematics, Lviv National University,
Universytetska 1, 79000 Lviv, Ukraine, and Institute of Mathematics,  University of Rzesz\'ow,
Rzesz\'ow, Poland} \email{topos@franko.lviv.ua,
mzar@litech.lviv.ua}

\date{\today}
\thanks{}
\begin{abstract} The objects of the Dranishnikov asymptotic category are proper metric spaces and the morphisms are asymptotically Lipschitz maps.
In this paper we  provide an example of an asymptotically  zero-dimensional space (in the sense of Gromov) whose space of
compact convex subsets of probability measures is not an absolute extensor in the asymptotic category in the sense of Dranishnikov.
\end{abstract}

\maketitle

\section{Introduction}

The notion of absolute extensor plays an important role in different branches of mathematics. In  asymptotic topology,
the absolute extensors are used in constructing the homotopy theory and the asymptotic dimension theory. Among the two
categories widely used in  asymptotic category, the {\it Dranishnikov} and the {\it Roe} categories (see the definition below),  it turns out that it is the
Dranishnikov category (the category of proper metric spaces and the asymptotically
Lipschitz maps) in which a richer extensor theory can be developed.

It was proved
in \cite{Za}  that in general, the space of
probability measures of a metric space is not an absolute extensor
for the Dranishnikov category. This provided a negative answer to a  question
formulated by Dranishnikov \cite[Problem 12]{d}, in connection with
existence of the homotopy extension theorem in this category. This leads to an open problem of searching functorial constructions that preserve  the class of absolute extensors in the asymptotic categories.

In the present
paper we deal with the hyperspaces of  compact convex subsets of probability measures. 
Note that these hyperspaces  play an important role in the decision theory, mathematical economics and finance, in particular, in the maximum (maxmin) expected utility  theory (cf. e.g. \cite{gs}).

In the case of compact metric spaces as well as in the case of compact spaces of weight $\omega_1$, the hyperspaces of compact convex subsets of probability measures  are known to be absolute extensors \cite{BRZ}. However, the extension properties of these hyperspaces in the asymptotic category  remained unknown. Our aim is to demonstrate that the example presented in \cite{Za} also works for the hyperspaces compact convex subsets of probability measures. Thus the main result of this paper
is
that the  spaces
mentioned above are not in general, asymptotic
extensors in the asymptotic category.

\section{Preliminaries}

\subsection{Asymptotic category}

Together with Roe's category of proper metric spaces and coarse
maps \cite{r}, the asymptotic category $\mathcal A$  introduced by
Dranishnikov \cite{d} turned out to be an important universe
for developing asymptotic topology.

A typical metric will be denoted by $d$. A map $f\colon X\to Y$ between metric
spaces is called $(\lambda,\varepsilon)$-{\it Lipschitz\/} for $\lambda>0$,
$\varepsilon\ge0$ if $d(f(x),f(x'))\le\lambda d(x,x')+\varepsilon$ for every
$x,x'\in X$. A map is called {\it asymptotically Lipschitz\/} if it is
$(\lambda,\varepsilon)$-Lipschitz for some $\lambda,\varepsilon>0$.
The $(1,0)$-Lipschitz maps are also called  {\it short}. The set of all short functions on a metric space $X$ is denoted by $\mathrm{LIP}(X)$.

A metric space is {\em proper} if every  closed ball in it is
compact. A map of metric spaces is (metrically) {\em proper} if
the preimages of the bounded sets are bounded. The objects of the
category $\mathcal A$ are the proper metric spaces and the
morphisms are the proper asymptotically Lipschitz maps.

A metric space $Y$ (not necessarily an object of $\mathcal A$) is
an {\em absolute extensor} (AE) for the category $\mathcal A$ if
for every proper asymptotically Lipschitz map $f\colon A\to Y$
defined on a closed subset of a proper metric space $X$ there
exists a proper asymptotically Lipschitz extension $\bar f\colon
X\to Y$ of $f$.

\subsection{Asymptotic dimension}

The notion of asymptotic dimension was
introduced by Gromov
\cite{g}.
Let $X$ be a metric space. A family $\mathcal C$ of subsets of $X$ is said to be
{\em uniformly bounded} if there exists $M>0$ such that $\diam A\le M$ for every
$A \in \mathcal C$. Given $D>0$, we say that a family $\mathcal C$
of subsets of $X$ is {\em $D$-disjoint} if
$\inf\{d(a,a')\mid a\in A,\ a'\in A'\}>D$ for every $A,A'\in \mathcal C$,
$A\neq A'$.

We say that the {\em asymptotic dimension} of $X$ is $\le n$
(written $\asd X\le n$) if for every $D>0$ there exists a  cover $\U$ of $X$
such that $\U=\U^0\cup\dots\cup\U^n$, where every family $\U^i$ is
$D$-discrete. If we require
in the  definition of the absolute extensor  that
$\asd X\le n$, then the definition of the {\em absolute extensor
in asymptotic dimension $n$} (briefly AE($n$)) is obtained.

It is easy to see that for a proper metric space $X$,  the
inequality $\asd X\le0$ is equivalent to the condition that for
every $C>0$ the diameters of the $C$-chains in $X$ (i.e. the
sequences $x_1,x_2,\dots,x_k$ with $d(x_i,x_{i+1})\le C$ for every
$i=1,2,\dots,k-1$) are bounded from above.

\subsection{Convex hyperspaces of  probability measures}

Let $P(X)$ denote the space of probability measures of compact supports on a metrizable space $X$. For any $x\in X$,  we denote the Dirac measure concentrated at $x$
by $\delta_x$.
If $d$ is a metric on $X$, we denote by $\hat d$ the Kantorovich metric generated by $d$,
$$\hat d(\mu,\nu)=\sup\left\{\left|\int\varphi d\mu -\int\varphi d\nu\right|\mid \varphi\in \mathrm{LIP}(X)\right\}$$
(cf. e.g. \cite{h}).

By $\mathrm{cc}P(X)$ we denote the set of all nonempty compact convex subsets in $P(X)$; as usual, a subset $A\subset P(X)$ is convex if $t\mu+(1-t)\nu\in A$, for all $\mu,\nu\in P(X)$ and $t\in[0,1]$. The set $\mathrm{cc}P(X)$ is endowed with the Hausdorff metric, which we shall
denote by $\hat d_H$:
$$\hat d_H(A,B)=\inf\{r>0\mid A\subset O_r(B),\ B\subset O_r(A)\}$$
(here $O_t(Y)$ stands for the $t$-neighborhood of $Y\subset P(X)$).
Note that, clearly, the map $x\mapsto\{\delta_x\}\colon X\to\mathrm{cc}P(X)$ is an isometric embedding.

Given a map $f\colon X\to Y$ of metric spaces, we define the map $P(f)\colon P(X)\to P(Y)$ as follows: $\int \varphi dP(f)(\mu)= \int \varphi f d\mu$.
The map $\mathrm{cc}P(f)\colon \mathrm{cc}P(X)\to\mathrm{cc}P(Y)$ is then defined by the formula:
$$\mathrm{cc}P(f)(A)=\{P(f)(\mu)\mid \mu\in A\}.$$ It can be easily seen that the map $\mathrm{cc}P(f)$ is short if such is $f$.

Let $b\colon P(\mathbb R^n)\to \mathbb R^n$ denote the barycenter map. Recall that this map assigns to every $\mu\in P(\mathbb R^n)$ the unique point $b(\mu)$ with the property that $L(b(\mu))=\int Ld\mu$, for every continuous linear functional $L$ on $\mathbb R^n$. Since $b$ is known to be continuous and linear, the image $b(A)$ of every $A\in\mathrm{cc}P(\mathbb R^n)$ is a compact convex subset of $\mathbb R^n$, i.e.,  an element of the space $ \mathrm{cc}(\mathbb R^n)$ of compact convex subsets in $\mathbb R^n$ endowed with the Hausdorff metric.

Let $p\colon \mathrm{cc}(\mathbb R^n)\to \mathbb R^n$ denote  the map defined by the condition:
$$y=\pi(A)\ \Leftrightarrow\ y\in A\text{ and }\|y\|=\inf\{\|z\|\mid z\in A \}.$$

The proof of the following statement uses simple geometric arguments and will therefore be omitted.
\begin{lemma}
The map $\pi$ is well-defined and short.
\end{lemma}

\section{The Example}\label{s:e}

Our example described  below is a modification of the second author's example
of a proper metric space whose space of probability measures is
not an AE (even AE(0)) in the asymptotic category \cite{z}. For the sake of
completeness we shall provide here the details of the construction.

For every $n$, the Euclidean space $\R^n$ can naturally be identified
with the subspace $\{(x_i)\mid x_i=0\text{ for all }j>n\}$ of the
space $\ell^2$.
We endow the subspace $X'=\bigcup_{n\in\N}\{n^2\}\times\R^n\subset \R\times\ell^2$
with the metric $$d((m,(x_i)), (n,(y_i)))=(|m-n|^2+\|(x_i)-(y_i)\|^2)^{1/2}.$$
Obviously, $X$ is a proper metric space. For every $n$ we denote by $p_n\colon
X'\to\R^n$ a map defined by the formula $p_n(m,(x_i))=(x_1,\dots,x_n)$. Clearly,
$p_n$ is a short map.

It was shown in \cite{l} (cf. Theorem 1.5 therein) that for any
$n\ge2$ there exists a metric space  $X_n$ which
contains  the
Euclidean space $\R^n$ as a metric subspace and such that there is
no $(\lambda,\varepsilon)$-Lipschitz retraction from $X_n$ onto
$\R^n$ with $\lambda<n^{1/4}$. In the sequel we shall need an
explicit construction of these spaces.
 Following \cite{l}, for every natural $k$ and natural $n\ge2$
we define graphs $G_{n,k}$ as follows: the set of vertices $V(G_{n,k})$ is the
union of $I(G_{n,k})$ and $T(G_{n,k})$, where
\begin{align*} I(G_{n,k})=&\{x=(x_1,\dots,x_n)\in\R^n\mid |x_i|=k\text{ for all
}i\},\\
T(G_{n,k})=&\{x=(x_1,\dots,x_n)\in\R^n\mid |x_i|=2k\text{ for all }i\};
\end{align*}
the set of edges $E(G_{n,k})$ is defined by the condition: $\{x,y\}\in
E(G_{n,k})$ if and only if $x,y\in V(G_{n,k})$ and either $\|x-y\|=2k$ or $y=2x$
(we suppose that the spaces $\R^n$ are endowed with the Euclidean metric).

The set  $V(G_{n,k})$ is equipped with the metric $d=d_{n,k}$,
\begin{align*} d(x,y)=&\inf\left\{\sum_{i=1}^l\|x_{i-1}-x_i\|_\infty\mid
(x=x_0,x_1,\dots,x_l=y)  \text{ is a path in }G_{n,k}\right\}
\end{align*}
(as usual, $\|x\|_\infty$ denotes the max-norm of $x\in\R^n$.)

Define spaces $X$ and $Y$ as follows:
 $$X=\bigcup_{n=2}^\infty \bigcup_{k=n}^\infty\{n^2\}\times T(G_{n^2,k^2}),\
 Y=\bigcup_{n=2}^\infty\bigcup_{k=n}^\infty \{n^2\}\times V(G_{n^2,k^2}) $$
where the metric on $X$ is inherited from $X'$ and  the metric on $Y$ is the
maximal metric that agrees with the already defined metric on $X$ and the metric
 $d_{n,k}$ on every
$V(G_{n,k})$. It easily follows from the construction that $X$ and $Y$ are
proper metric spaces, i.e. objects of the category $\mathcal A$.

We are going to show that $\asd Y=0$ (and consequently $\asd
X=0$). Let $C>0$ and suppose that $y_1,\dots,y_m$ is a $C$-chain
in $Y$. Denote by $k$ the minimal natural number such that
$C<(k+1)^2$. If
$$\{y_1,\dots,y_n\}\subset \bigcup_{j=2}^k\bigcup_{l=j}^k
\{j^2\}\times V(G_{j^2,l^2})$$ then $\diam\{y_1,\dots,y_n\}\le
\sqrt{(k^2)+(3k^2)^2}\le \sqrt{10C}$.
 Otherwise $$\{y_1,\dots,y_n\}\cap Y\setminus \bigcup_{j=2}^k\bigcup_{l=j}^k
\{j^2\}\times V(G_{j^2,l^2})\neq\emptyset$$ and
$\{y_1,\dots,y_n\}$ is a singleton.

It was proved in \cite{l} that  the following holds for the spaces
$$X_{n^2}= \R^{n^2}\cup\bigcup_{k=n}^\infty \{n^2\}\times
V(G_{n^2,k^2})$$ endowed with the maximal metric which
agrees with
the initial metric on $\R^{n^2}$ and the  metric
on $\bigcup_{k=n}^\infty \{n^2\}\times V(G_{n^2,k^2})$ inherited from $Y$
(note that
these two metrics coincide on the intersection of their domains):
there is no $(\lambda,\e)$-retraction of $X_{n^2}$ to $\R^{n^2}$
with $\lambda<\sqrt{n}$.

Now, let $f\colon X\to \mathrm{cc}P(X)$ be the map that sends $x\in X$ to
$\{\delta_x\}\in \mathrm{cc}P(X)$. The map $f$ is an isometric embedding and we
are going to show that there is no asymptotically Lipschitz
extension of $f$ onto the whole space $Y$. Assume the contrary and let $\bar
f\colon Y\to \mathrm{cc}P(X)$ be such an extension. We regard $\bar f$ as a
map into $F(X')\supset \mathrm{cc}P(X)$. Then there exist $\lambda>0$ and
$\varepsilon>0$ such that $$d(\bar f(x),\bar f(x'))\le\lambda
d(x,x')+\varepsilon$$ for all $x,x'\in Y$.

Let $n>\lambda^2$. Since the maps $\mathrm{cc}P(p_{n^2})$,  $b\colon
P(\R^{n^2})\to\R^{n^2}$ and $\pi$ are short, we conclude that the map $$x\mapsto \pi(\{b(\mu)\mid \mu\in \mathrm{cc} P(p_{n^2})(\bar f(x))\})\colon X_{n^2}\to\R^{n^2}$$ is a
$(\lambda,\varepsilon)$-Lipschitz retraction from $X_{n^2}$ onto
$\R^{n^2}$, which  contradicts to the choice of $\lambda$. This
demonstrates that the space $\mathrm{cc}P(X)$ is not an AE(0) for the asymptotic category $\mathcal A$.

\section{Epilogue}
We conjecture that the spaces of capacities (non-additive measures; cf. e.g. \cite{nz,zhou}) are always absolute extensors in the asymptotic category $\mathcal A$. Note that the scheme of our proof of the main result of this paper does not work for non-additive situation, because one does not have the ``barycenter map'' in this case.

\section*{Acknowledgements}

This research was supported by the
Slovenian Research Agency grants P1-0292-0101, J1-9643-0101 and
J1-2057-0101. We thank the referee for comments and suggestions.

\end{document}